\documentstyle{article}
\begin{document}
\begin{large}

\centerline{\bf Dopovidi of the Academy of Sciences of the Ukraine, Ser. A}
\centerline{\bf (1987), No. 10, p. 9-12.}
\centerline{\bf M.I.OSTROVSKII}
\centerline{$W^{*}$-{\bf DERIVED SETS OF TRANSFINITE ORDER OF SUBSPACES}}
\centerline{\bf OF DUAL BANACH SPACES}

Let $X$ be a separable Banach space and $X^*$  be  its  dual  space.
If $\Gamma $ is a linear (not necessarily closed) subspace of $X^{*}$, then the
$w^*${\it -derived set} of $\Gamma $, denoted $\Gamma _{(1)}$, is
the set of all limits of $w^{*}$-convergent sequences
in $\Gamma $. Generally speaking, $\Gamma _{(1)}$ may be non-closed
in $w^*$-topology [1]. This was a reason to S.Banach to  introduce  [2,
p.~213] derived sets of other orders, including  transfinite  ones.
Their inductive definition is the following. The $w^{*}$-{\it derived  set  of
order} $\alpha $ is the set
$$
\Gamma _{(\alpha )}=\bigcup _{\beta <\alpha }((\Gamma _{(\beta )})_{(1)}.
$$
We have $\Gamma \subset \Gamma _{(\beta )}\subset
\Gamma _{(\alpha )}$  for $\beta <\alpha $  and
if $\Gamma _{(\alpha )}=\Gamma _{(\alpha +1)}$  then  all
subsequent derived sets are equal to each other.
In [2, p.~209--213; 3--7] it was proved that for  some  separable
Banach spaces there exist subspaces in their  duals  for  which  the
chains of different $w^{*}$-derived sets is long in  certain  sense  (see
Commentary). We shall also consider only separable spaces.

Let us recall that Banach space $X$ is called {\it quasireflexive}   if
$\dim (X^{**}/\pi (X))<\infty $, where $\pi :X\to X^{**}$ is
the canonical embedding.

{\bf THEOREM}. Let $X$ be a nonquasireflexive separable  Banach  space.
Then for every countable ordinal $\alpha $ there is a
linear subspace $\Gamma \subset X^{*}$
for which $\Gamma _{(\alpha )}\neq \Gamma _{(\alpha +1)}=X^*$.

Remarks. Let $X$ be a separable  Banach  space.  Then  for  every
subspace $\Gamma \subset X^{*}$ there exists a
countable ordinal $\alpha $ for which
$\Gamma _{(\alpha )}=\Gamma _{(\alpha +1)}$
[8, p.~50]. In [6] it is proved that the length of the chain of
different $w^{*}$-derived sets cannot be equal to a limit ordinal. From
[2, p.~213] and well-known properties of quasireflexive spaces (see
[8, Chapter 4]) it follows that if $X$ is a quasireflexive separable
Banach space then for
every subspace $\Gamma \subset X^{*}$
we have $\Gamma _{(1)}=\Gamma _{(2)}$. Hence
our result is the best possible for separable Banach spaces.

{\bf Lemma 1}. Let $X$ be a separable Banach space, $Y$ be its subspace,
$\xi :Y\to X$ be the operator of
the identical embedding, and $\Gamma $ be a
subspace of $Y^{*}$. Then for every ordinal $\alpha $ we have
$(\xi ^{*})^{-1}(\Gamma _{(\alpha )})=((\xi ^{*})^{-1}\Gamma )_{(\alpha )}$.

Proof. The
inclusion
$(\xi ^{*})^{-1}(\Gamma _{(\alpha )})\subset
((\xi ^{*})^{-1}\Gamma )_{(\alpha )}$
follows by the
$w^{*}$-continuity of the operator $\xi ^{*}$.
In  order  to  prove  the  inverse
inclusion it is sufficient to show   that   for   every   sequence
$\{f_i\}^{\infty }_{i=1}\subset Y^{*}$ with $w^*-\lim f_i=f$ and
any $g\in (\xi ^{*})^{-1}(\{f\})$ there exist
$g_{i}\in (\xi ^{*})^{-1}(\{f_i\})$ such that $w^{*}-\lim g_{i}=g$.
Let $\{e_{i}\}^{\infty }_{i=1}$ be a sequence for which
$e_{i}\in (\xi ^{*})^{-1}(\{f_{i}\})$ and $||e_{i}||=
||f_{i}||$. Then the
sequence $\{e_{i}-g\}^{\infty }_{i=1}$ is
a bounded sequence of $X^{*}$ and
all its limit points
are in $(\xi ^{*})^{-1}(\{0\})$. Weak$^{*}$ topology is metrizable
on bounded subsets of the dual of separable Banach space.  Therefore
there exists a sequence
$\{h_{i}\}^{\infty }_{i=1}\subset (\xi ^{*})^{-1}(\{0\})$
for which $w^{*}-\lim (e_{i}-g-h_{i})=0$. It
is clear that the
vectors $g_{i}=e_{i}-h_{i} (i\in {\bf N})$ forms the desired sequence.

Proof of the theorem. We need the  following  fact [9]:  every
nonquasireflexive Banach space $X$ contains  a  bounded  away  from  0
basic sequence $\{z_{n}\}^{\infty }_{n=0}$ for which $||z_n||\le 1$
and the set

$$
\{||\sum^{k}_{i=j}z_{i(i+1)/2+j}||\}^{\infty }_{j=0},^{\infty }_{k=j}
$$
is bounded. Let us denote by $Z$ the closure of the linear span of the
sequence $\{z_{n}\}^{\infty }_{n=0}\subset X$. In order
to prove the theorem it is sufficient by
Lemma 1 to find a
subspace $\Gamma \subset Z^{*}$ for
which $\Gamma _{(\alpha )}\neq \Gamma _{(\alpha +1)}=Z^{*}$.

Let us introduce some notations. We will write $z^{j}_{i}$ for
$z_{(j+i-1)(j+i)/2+j}$, biorthogonal
functionals of the system $z_{n} (z^{j}_{i})$
will be denoted by $\tilde{z}_{n} (\tilde{z}^{j}_{i})$.
By abovementioned result from [9] we
have
$$\sup_{j,m}||\sum^{m}_{i=1}z^{j}_{i}||=M_{1}<\infty.$$
Therefore for every $j=0,1,2,\ldots $ the sequence
$\{\sum^{m}_{i=1}z^{j}_{i}\}^{\infty }_{m=1}$ has at
least one $w^{*}$-limit point in $Z^{**}$. Let us for every
$j=0,1,2\ldots $ choose one of such limit points and denote it by $f_{j}$. It
is clear that $||f_{j}||\le M_{1}$.

{\bf Lemma 2.} For every vector $g_{0}\in Z^{**}$ of
the form $af_{j}+z^{r}_{s} (a>0, r\neq j)$,
every countable ordinal $\alpha $ and
every infinite subset $A\subset {\bf N}$ such that
$j,r\not\in A$ there exists a
countable subset $\Omega (g_{0},\alpha ,A)\subset Z^{**}$ such that

1) The set $K(g_{0},\alpha ,A)$ defined by
$K(g_{0},\alpha ,A)=(\cap \{\hbox{ker}h:\ h\in \Omega (g_{0},\alpha ,A)\})$
satisfies the condition: $(K(g_{0},\alpha ,A))_{(\alpha )}\subset \ker  g_{0}.$

2) All vectors $h\in \Omega (g_{0},\alpha ,A)$ are of  the
form $h=a(h)f_{j(h)}+z^{r(h)}_{s(h)}$
with $j(h), r(h)\in A\cup \{j,\ r\}, a(h)>0,$ and if we have
$j(h)=r$ or $r(h)=r$ then $h=g_{0}$.

3) Every finite linear combination of vectors
$\{\tilde{z}_{n}\}$ which is  of
the form
$$
b\tilde{z}^{r}_{s}+u; u\in
\hbox{lin}(\{\tilde{z}^{t}_{k}\}^{\infty }_{k=1,t\in A\cup \{j\}})
\eqno{(1)}$$
 must be in $(Q(b,g_{0},\alpha ,A)_{(\alpha )}$
(where $Q(b,g_{0},\alpha ,A)$ is the  set  of  all
linear combinations of the type (1) which are in
$K(g_{0},\alpha ,A))$ if it is
in $\ker  g_{0}.$
(Here we need to remark that the definition of the $w^{*}$-derived
 sets for subsets which are not subspaces is the same.)

At first we will finish the proof of the theorem with the help
of lemma 2.
>From 1) it follows that $(K(g_{0},\alpha ,A))_{(\alpha )}\neq Z^{*}$.
Let us show that 3)
implies $(K(g_{0},\alpha ,A))_{(\alpha +1)}=Z^{*}$. Let us recall that
$\{z_{n}\}^{\infty }_{n=0}$ is a basis
 in $Z$, therefore $\{\tilde{z}_{n}\}^{\infty }_{n=0}$ is a
$w^{*}$-Schauder basis in $Z^{*}$ [10, p. 155],
 i.e. every vector $\tilde{z}\in Z^{*}$ can be represented as
$\tilde{z}=w^{*}-\lim_{n\to\infty}\sum^n_{k=1}a_{k}\tilde{z}_{k}$, where
$a_{k}=\tilde{z}(z_{k}).$

It  is  clear  that  vectors
$\tilde{y}_{n}=\sum^{n}_{k=0}a_{k}\tilde{z}_{k}-g_{0}
(\sum^{n}_{k=0}a_{k}\tilde{z}_{k})\tilde{z}^{j}_{n}/a$
can   be represented in form $\tilde{y}^{1}_{n}+\tilde{y}^{2}_{n}$
where $\tilde{y}^{1}_{n}$ is of type (1) and is in $\ker  g_{0}$,
and $\tilde{y}^{2}_{n}$ is a finite linear combination of the vectors
$\{\tilde{z}_{n}\}^{\infty }_{n=0}\backslash
(\{\tilde{z}^{t}_{k}\}^{\infty }_{k=1,t\in A\cup
\{j\}}\cup \{\tilde{z}^{r}_{s}\}),$
and, hence, $\tilde{y}^{2}_{n}\in K(g_{0},\alpha ,A).$ By 3) we
have $\tilde{y}^{1}_{n}\in (K(g_{0},\alpha ,A))_{(\alpha )}.$
Therefore $\tilde{y}_{n}\in (K(g_{0},\alpha ,A))_{(\alpha )}.$ It is clear
 also that $w^{*}-\lim \tilde{y}_{n}=\tilde{z}.$ Thus the proof of the
theorem is complete.

Proof of Lemma 2. A. Let us suppose that the assertion of Lemma
2 is true for ordinal $\alpha $. Let us show that it is true for
$\alpha +1.$ Let us
represent the set $A$  as  a  countable  union  of pairwise disjoint
infinite subsets, $A=\cup ^{\infty }_{k=0}A_{k}$.
Let $\varepsilon _{i}>0\ (i\in {\bf N})$ be such that
$\sum^{\infty }_{i=1}\varepsilon _{i}<\infty $.
 Let us introduce the sequence $\{g_{n}\}^{\infty }_{n=1}\subset Z^{**}$ by
$g_{n}=\varepsilon _{n}f_{p(n)}+z^{j}_{n}$, where
$p:{\bf N}\to A_{0}$ is some one-to-one mapping.
Let $\Omega (g_{n},\alpha ,A_{n})$ be the  sets  whose
existence follows from our assumption. Let us show that we can let
$\Omega (g_{0},\alpha +1,A)=(\cup ^{\infty }_{n=1}\Omega 
(g_{n},\alpha ,A_{n}))\cup \{g_{0}\}.$ It is clear that 2) is satisfied.

 Let us show that 1) is satisfied.  In  order  to  do  this  we  will
using  transfinite  induction  show   that   for 
$\beta \le \alpha +1$   we   have
$(K(g_{0},\alpha +1,A))_{(\beta )}\subset \ker  g_{0}$. 
Let us  suppose  that  we  prove  this  for
certain $\beta \le \alpha $.   Let   us   derive   from   this   
assumption   that
$(K(g_{0},\alpha +1,A))_{(\beta +1)}\subset \ker  g_{0}$. 
In order to do this we  must  prove  that
for any $w^{*}$-convergent sequence $\{\tilde{y}_{i}\}^{\infty }_{i=1}\subset 
(K(g_{0},\alpha +1,A))_{(\beta )}$ its limit
 $\tilde{y}=w^{*}-\lim\tilde{y}_{i}$ belongs to $\ker  g_{0}$. 
We have already noted that $\{\tilde{z}_{n}\}$ is
$w^{*}$-basis in $Z^{*}$. Let us estimate coefficients 
$\{\alpha (i)^{j}_{n}\}^{\infty }_{n=1}$ which are
staying near vectors $\{\tilde{z}^{j}_{n}\}^{\infty }_{n=1}$ 
in the corresponding $w^{*}$-decomposition
of $\tilde{y}_{i}$. It  is  clear  that $\sup||\tilde{y}_{i}||=M_{2}<\infty $.  
>From $\beta \le \alpha $  and  induction
hypothesis we obtain $g_{n}(\tilde{y}_{i})=0.$ We have 
$g_{n}(\tilde{y}_{i})=\alpha (i)^{j}_{n} +\varepsilon _{n}f_{p(n)}
(\tilde{y}_{i})$,
 hence
$$
| \alpha (i)^{j}_{n}| \le \varepsilon _{n}M_{1}M_{2} (n\in {\bf N}).
\eqno{(2)}$$
Therefore vector $\tilde{y}_{i}$ can be represented in the form 
$\tilde{y}_{i}=u_{i}+v_{i}$
in such a way that $w^{*}$-decomposition of $u_{i}$ doesn`t involve vectors
$\tilde{z}^{r}_{s}$ and $\{\tilde{z}^{j}_{n}\}^{\infty }_{n=1}$
and $v_{i}=a_{i}\tilde{z}^{r}_{s}+\sum^{\infty }_{n=1}
\alpha (i)^{j}_{n}\tilde{z}^{j}_{n}.$ By (2) this  series  is
strongly  convergent.  It  is  clear  that $w^{*}$-convergence   implies
coordinatewise convergence, therefore it follows  from  (2)  that $\tilde{y}$
also can be represented in the form $\tilde{y}=u+v$ in the same manner. So we
have $v=w^{*}-\lim v_{i}$ and $u=w^{*}-\lim u_{i}.$ From $u_{i}(\sum^{m}_{k=1}
\gamma _{k}z^{j}_{k}+\delta z^{r}_{s})=0\ (i\in {\bf N})$ and
analogous equality for $u$ which holds for all 
$m, \{\gamma _{k}\}^{m}_{k=1}$  and $\delta $  it
follows that $g_{0}(u_{i})=0,\  g_{0}(u)=0.$ 
Since $\tilde{y}_{i}\in \ker  g_{0}$ we have $v_{i}\in \ker  g_{0}$.

Thus it is sufficient to prove that $v\in \ker  g_{0}$. 
This assertion follows
from $v=w^{*}-\lim v_{i};\ v_{i}\in \ker g_{0}$ and the following two facts:

a)  For   every $i\in {\bf N}$   we   have 
$v_{i}\in V:=\{\tilde{z}\in Z^{*}:\tilde{z}=
\sum^{\infty }_{k=1}c_{k}\tilde{z}^{j}_{k}+d\tilde{z}^{r}_{s},
| c_{k}| \le \varepsilon _{k}M_{1}M_{2},\ |d|\le M_{2}\}$ 
and the set $V$ is compact in the strong topology.

b) Initial topology of compact space coincides with  any
other Hausdorff topology which is weaker than it [11, p.~249].

Therefore we proved that if $\beta \le \alpha $ and 
$(K(g_{0},\alpha +1,A))_{(\beta )}\subset \ker  g_{0}$ then
$$
(K(g_{0},\alpha +1,A))_{(\beta +1)}\subset \ker  g_{0}.
$$
The case of limit ordinal is evident.

The  proof  of  3)  follows  immediately  from  the   induction
hypothesis  and  the  following  assertions,  which  can  be  easily
verified:

(i)  any  vector  of  type  (1)  can  be  represented  in  form
$b\tilde{z}^{r}_{s}+\sum^{m}_{k=1}a_{k}\tilde{z}^{j}_{k}+u_{k}$, where
$u_{k}\in $lin $(\{\tilde{z}^{t}_{l}\}^{\infty }_{l=1,t\in A_k\cup \{p(k)\}})$.

(ii) if $M$ contains all vectors of the type (1) from $\ker  g_{0}$ then
 $M_{(1)}$ contains all vectors of the type (1).

(iii)  if  collection $\{a_{k}\}^{m}_{k=1}\subset {\bf R}$ 
and $b\in {\bf R}$   are   such   that
$g_{0}(b\tilde{z}^{r}_{s}+\sum^{m}_{k=1}a_{k}\tilde{z}^{j}_{k})=0$ then
we have $Q(b,g_{0},\alpha +1,A)\supset
b\tilde{z}^{r}_{s}+\sum^{m}_{k=1}Q(a_{k},g_{k},\alpha ,A_{k})$.

(iv) for any finite collection of sets $M_{i}\subset Z^{*}, i=1,2,\ldots
,m$, and any ordinal $\alpha $ we have 
$(\sum^{m}_{i=1}M_{i})_{(\alpha )}\supset \sum^{m}_{i=1}(M_{i})_{(\alpha )}$.

B. The case in which $\alpha $ is a limit ordinal. By [11,  p.~72]  we
can find an increasing sequence $\{\alpha _{i}\}^{\infty }_{i=1}$ of ordinals 
for which $\alpha $ is a limit.

Let $\{A_{n}\}^{\infty }_{n=0}$ and $\{g_{n}\}^{\infty }_{n=1}$ 
are the same  as  in  $A$.  By  induction
hypothesis we can find subsets 
$\Omega (g_{n},\alpha _{n},A_{n})\subset Z^{**}$. 
The fact that $(\cup ^{\infty }_{n=1}\Omega (g_{n},\alpha _{n},A_{n}))\cup 
\{g_{0}\}$ is the desired set can be proved in the same
way as in $A$ with unique distinction that inequalities (2)  are valid
not for all $n$ but for all exept finite number of them.

Received January 19.1987.

Address: Institute for low Temperature Physics and
\par
Engineering
\par
Ukrainian Academy of Sciences
\par
47 Lenin avenue
\par
Kharkov 310164 Ukraine
\par
USSR.
\par
REFERENCES

1.  Mazurkiewicz  S.  Sur la  derivee  faible d`un   ensemble   de
fonctionnelles lineaires, Stud. Math. 2. (1930), 68-71.

2.  Banach  S.  Theorie   des   operations   lineaires.   Monografje
Matematyczne. no. 1 (Warszawa, 1932).

3. McGehee O.C. A proof of a statement of  Banach  about  the  weak$^{*}$
topology, Mich. Math. J. 15 (1968), no. 2, 135-140.

4. Sarason D. On the order of a simply connected domain, Mich. Math.
J. 15 (1968), no. 2, 129-133.

5. Sarason D. A remark on the weak-star topology of $l^{\infty }$, Stud.  Math.
30 (1968), no. 3, 355-359.

6. Godun B.V. On weak$^{*}$ derived sets of transfinite order of sets  of
linear functionals. Sib. Mat. Zh. 18  (1977),  no.  6,  1289-1295
(Russian).

7. Godun B.V. On weak$^{*}$ derived sets of sets of  linear  functionals.
Mat. Zametki 23 (1978), no. 4, 607-616 (Russian).

8. Petunin Yu.I., Plichko A.N.  The  theory  of  characteristics  of
subspaces and its applications (Kiev, Vyshcha Shkola, 1980).

9. Davis W.J., Johnson W.B. Basic sequences and norming subspaces in
non-quasi-reflexive   Banach   spaces,   Isr.   J.    Math.    14
(1973),353-367.

10. Singer I. Bases in Banach spaces, I (Berlin, Springer, 1970).

11. Aleksandrov P.S. Introduction to  the  set  theory  and  general
topology (Moscow, Nauka, 1977).

\centerline {\bf COMMENTARY BY THE AUTHOR (July, 1991).}

I  would  like  to  say  several  words  about   history   and
applications of ``long chains'' of $w^{*}$-derived sets.

S.Mazurkiewicz [1] was the first
who found a subspace $\Gamma \subset (c_{0})^{*}$
for  which $\Gamma _{(1)}\neq \Gamma _{(2)}$.  S.~Banach  [2, p.~209]
developed   this
construction and for any $n\in {\bf N}$ found a
subspace $\Gamma \subset (c_{0})^{*}$ for which
$\Gamma _{(n)}\neq \Gamma _{(n+1)}$ and stated
analogous result for any countable ordinal $\alpha $.
But S.Banach's proof of this result never appeared.  This  statement
has been proved by O.C.McGehee [3]. (It  is  interesting  to  remark
that O.C.McGehee used results from Fourier analysis.) At almost  the
same time D.Sarason (using the complex analysis) obtained  analogous
results for $l_{1}$ and $L_{1}/H_{1}$ [4, 5]. In [6] B.V.Godun
observed that  the
length of the chain of different $w^{*}$-derived sets must  be  countable
and cannot be equal to a limit ordinal. In [7] B.V.Godun proved that
for any nonquasireflexive separable
Banach space $X$ and any $n\in {\bf N}$ there
exists a total subspace $\Gamma \subset X^{*}$
for which $\Gamma _{(n)}\neq X^{*}$.

The problem of existence of a total subspace $\Gamma \subset X^{*}$ for which
$\Gamma _{(n)}\neq X^{*}$ for all $n\in {\bf N}$ is
turned out to be important in the  theory  of
topological vector spaces [12].  This  problem  has  been  posed  by
V.B.Moscatelli at the Ninth  Seminar  (Poland  -  GDR)  on  Operator
Ideals and Geometry of Banach spaces (Georgenthal, April, 1986) [13,
Problem 17]. This problem is solved by the result of  present  paper
(At the time of writing it the author had not any information  about
[12]  and  [13])  Independently  and  almost  at   the   same   time
V.B.Moscatelli [14] proved that for any separable  nonquasireflexive
Banach space $X$ there exists a subspace $\Gamma \subset X^{*}$
for which $\Gamma _{(\omega )}\neq X^{*}$ (where
$\omega $ is the first infinite ordinal). In  [15]  V.B.Moscatelli  obtained
this result in more explicit form. The reader may consult  [12,  15,
16]  for  applications  of  the  ``long  chains''  in  the  theory  of
topological vector spaces. In [17] ``long chains'' are used to answer
Kalton's question on universal biorthogonal systems. There are  also
certain applications of ``long chains'' in  the  theory  of  ill-posed
problems, or, more precisely, in the problems of regularizability of
inverses of injective continuous linear operators (see [8, 18, 19]).

Some  connections  of  ``long  chains''  with  harmonic  analysis  are
discussed in [20].

12. Dierolf S., Moscatelli V.B. A note  on  quojections,  Funct.  et
Approxim. (1987), n. 17, 131-138.

13. Open Problems, Presented at the Ninth Seminar (Poland - GDR)  on
Operator Ideals  and  Geometry  of  Banach  Spaces,  Georgenthal,
April, 1986, Forschung. Friedrich - Schiller - Universitat, Jena,
N/87/28, 1987.

14. Moscatelli V.B. On strongly non-norming subspaces, Note  Mat.  7
(1987), 311-314.

15.   Moscatelli   V.B.   Strongly    nonnorming    subspaces    and
prequojections, Studia Math. 95 (1990), 249-254.

16. Metafune G., Moscatelli V.B. Quojections and prequojections, in:
Advances in the Theory of Frechet  Spaces,  ed.  by  T.Terzioglu,
Kluwer Academic Publishers, Dordrecht, 1989, 235-254.

17. Plicko A.N. On bounded biorthogonal  systems  in  some  function
spaces, Studia Math. 84 (1986), 25--37.

18. Ostrovskii M.I. Pairs of regularizable inverse linear  operators
with nonregularizable superposition, J.Soviet  Math.  52  (1990),
no. 5, 3403--3410.

19. Ostrovskii M.I.,  Regularizability   of
superpositions of inverse  linear  operators,  Teor.  Funktsii,
Funktsional. Anal i Prilozhen. 55 (1991), 96--100 (in Russian). Engl.
transl.: J. Soviet Math. 59 (1992), no. 1, 652--655.

20.  Katznelson  Y.,  McGehee  O.C.  Some  sets   obeying   harmonic
synthesis, Israel J. Math. 23 (1976), no. 1, 88--93,

\end{large}
\end{document}